\documentclass[12pt]{article}
\usepackage{amssymb}
\usepackage{amsfonts}
\usepackage{latexsym}
\usepackage[dvips]{graphics}

\newtheorem{thm}{Theorem}
\newtheorem{lem}{Lemma}

\newtheorem{cor}{Corollary}

\newtheorem{rem}{Remark}
\newtheorem{prob}{Problem}

\large\normalsize

\begin{document}

\begin{center}
%%%%%%%%% title %%%%%%%%%%
{\Large\bf The determinants of $q$-distance matrices of trees and
two quantities relating to permutations}
\\[30pt]
%%%%%%%%%% authors %%%%%%%%%%
{Weigen\ Yan$^{\rm a,b}$ \footnote{This work is supported by
FMSTF(2004J024) and NSFF(E0540007).} \quad and \quad Yeong-Nan\
Yeh$^{\rm b}$ \footnote{Partially supported by
NSC95-2115-M001-009.
\newline\hspace*{5mm}{\it Email address:} weigenyan@263.net (W. G.
Yan), mayeh@math.sinica.edu.tw (Y. N. Yeh).}}
%\footnote {Email address: weigenyan@263.net (W.
%Yan), mayeh@math.sinica.edu.tw (Y. N. Yeh)}
\\[10pt]
%%%%%%%%%% affiliation %%%%%%%%%%
\footnotesize { $^{\rm a}$School of Sciences, Jimei University,
Xiamen 361021, China
\\[7pt]
$^{\rm b}$Institute of Mathematics, Academia Sinica, Taipei 11529.
Taiwan}
%\\[30pt]
\end{center}
%%%%%%%%%% abstract %%%%%%%%%
\begin{abstract}
In this paper we prove that two quantities relating to the length
of permutations defined on trees are independent of  the
structures of trees. We also find that these results are closely
related to the results obtained by Graham and Pollak (Bell System
Tech. J. 50(1971) 2495--2519) and by Bapat, Kirkland, and Neumann
(Linear Alg. Appl. 401(2005) 193--209).
\end{abstract}

%%%%%%%%%%%%%%%%%%%%%%%%%%%%%%%%%%%%%%%%%%
%%%%%%%%%%%%% Section 1
%%%%%%%%%%%%%%%%%%%%%%%%%%%%%%%%%%%%%%%%%%
\section{Introduction}
\hspace*{\parindent} Let $[n]$ denote the set $\{1,2,\ldots,n\}$
and let $\mathcal S_n$ be the set of permutations of $[n]$.
Partition $\mathcal S_n$ into $\mathcal S_n=\mathcal E_n\cup
\mathcal O_n$, where $\mathcal E_n$ (resp. $\mathcal O_n$) is the
set of even (resp. odd) permutations in $\mathcal S_n$. It is well
known that $|\mathcal E_n|=|\mathcal O_n|$. Let $\sigma$ and $\pi$
be two elements of $\mathcal S_n$. Diaconis and Graham \cite{DG77}
defined a metric called Spearman's measure of disarray on the set
$\mathcal S_n$ as follows:
$$D(\sigma,\pi)=\sum_{i=1}^n|\sigma(i)-\pi(i)|.$$
They derived the mean, variance, and limiting normality of
$D(\sigma,\pi)$ when $\sigma$ and $\pi$ are chosen independently
and uniformly from $\mathcal S_n$. In particular, the authors in
\cite{DG77} characterized those permutations $\sigma\in \mathcal
S_n$ for which $D(\sigma)=:D(1,\sigma)$ takes on its maximum
value. Some related work appears in \cite{HY97,Ster90}. The {\sl
length} $|\sigma|$ of a permutation $\sigma$ is defined to be
$D(1,\sigma)$, that is,
$|\sigma|=\sum\limits_{i=1}^n|i-\sigma(i)|$.
%Define the length of
%a permutation $\sigma$, denoted by $|\sigma|$, as $D(1,\sigma)$,
%\emph {i.e.}, $|\sigma|=\sum\limits_{i=1}^n|i-\sigma(i)|$.
For an arbitrary nonnegative integer $k$, let
$$\mathcal A_{n,k}=\{\sigma\in \mathcal S_n |\ |\sigma|=k\},$$
$$N_{n,k}=\sum_{\sigma \in \mathcal A_{n,k}}\mathrm{sgn}(\sigma)=
|\mathcal A_{n,k}\cap \mathcal E_n|-|\mathcal A_{n,k}\cap \mathcal
O_n|.$$ Furthermore, we define $\phi_{\sigma,k}=0$ if $\sigma$ has
at least one fixed point, otherwise, let $\phi_{\sigma,k}$ be the
number of nonnegative integer solutions of the equation
$x_1+x_2+\ldots+x_n=k$ which satisfy $0\leq x_i<|i-\sigma(i)|$ for
$1\leq i\leq n$. Let
$$M_{n,k}=\sum_{\sigma \in \mathcal S_{n}}\mathrm{sgn}(\sigma)\phi_{\sigma,k}. \eqno{(1)}$$
It is natural to pose the following problem:
%%%%%%%%%%%%%%%%%%%%%%%%%%%%%%%%%%%%%%%%%%
%%%%%%%%%%%%% Problem 1.1
%%%%%%%%%%%%%%%%%%%%%%%%%%%%%%%%%%%%%%%%%%
\begin{prob}
%For an arbitrary nonnegative integer $k$,
Find closed expressions for $N_{n,k}$ and $M_{n,k}$.
\end{prob}

We may generalize the concept of the length of a permutation
defined in Problem 1.1 as follows. Let $T$ be a weighted tree with
the vertex set $V(T)=\{v_1,v_2,\ldots,v_n\}$. For two vertices $u$
and $v$ in $T$, there exists a unique path $u=v_{i_1}$-$v_{i_2}$-\
$\ldots$\ -$v_{i_l}$-$v_{i_{l+1}}=v$ from $u$ to $v$ in $T$.
Define the {\sl distance} $d(u,v)$ between $u$ and $v$ as zero if
$u=v$, otherwise, let $d(u,v)$ be the sum $x_1+x_2+\ldots+x_l$,
where $x_k$ is the weight of edge $v_{i_k}v_{i_{k+1}}$ for
$k=1,2,\ldots,l$. Let $T$ be a simple tree (\emph{i.e.}, the
weight of each edge equals one) and let $\sigma\in \mathcal S_n$.
The {\sl length} $|\sigma_T|$ of $\sigma$ on $T$ is defined as the
sum of all $d(v_i,v_{\sigma(v_i)})$, that is,
$|\sigma_T|=\sum\limits_{i=1}^n d(v_i,v_{\sigma(i)})$. Let
$$\mathcal A_{n,k}(T)=\{\sigma\in \mathcal S_n |\ |\sigma_T|=k\},$$
$$N_{n,k}(T)=\sum_{\sigma \in \mathcal A_{n,k}(T)}\mathrm{sgn}(\sigma)=
|\mathcal A_{n,k}(T)\cap \mathcal E_n|-|\mathcal A_{n,k}(T)\cap
\mathcal O_n|.$$ Furthermore, we define $\phi_{\sigma,k}(T)=0$ if
$\sigma$ has at least one fixed point, otherwise, let
$\phi_{\sigma,k}(T)$ be the number of nonnegative integer
solutions of the equation $x_1+x_2+\ldots+x_n=k$ which satisfy
$0\leq x_i<d(v_i,v_{\sigma(i)})$ for $1\leq i\leq n$. Let
$$M_{n,k}(T)=\sum_{\sigma \in \mathcal S_{n}}\mathrm{sgn}(\sigma)\phi_{\sigma,k}(T).\eqno{(2)}$$
A more general problem than Problem 1.1 is the following:
%A more general problem than Problem 1.1 is the following:
%%%%%%%%%%%%%%%%%%%%%%%%%%%%%%%%%%%%%%%%%%
%%%%%%%%%%%%% Problem 1.2
%%%%%%%%%%%%%%%%%%%%%%%%%%%%%%%%%%%%%%%%%%
\begin{prob} Let $T$ be a simple tree with vertex set
$\{v_1,v_2,\ldots,v_n\}$.
%For an arbitrary nonnegative integer $k$,
Find closed expressions for $N_{n,k}(T)$ and $M_{n,k}(T)$.
\end{prob}
%%%%%%%%%%%%%%%%%%%%%%%%%%%%%%%%%%%%%%%%%%
%%%%%%%%%%%%% Remark 1.1
%%%%%%%%%%%%%%%%%%%%%%%%%%%%%%%%%%%%%%%%%%
\begin{rem}
If we take $T=P_n$ (where $P_n$ is a simple path with vertex
set $\{v_1,v_2,\ldots,v_n\}$ and edge set $\{(v_i,v_{i+1})| 1\leq
i\leq n-1\}$) in Problem 1.2, then Problem 1.1 is a special case
of Problem 1.2. That is, $N_{n,k}=N_{n,k}(P_n)$ and
$M_{n,k}=M_{n,k}(P_n)$.
\end{rem}

The distance matrix $D(T)$ of the weighted tree $T$ is an $n\times
n$ matrix with its $(i,j)$-entry equal to the distance between
vertices $v_i$ and $v_j$. If $T$ is a simple tree, Graham and
Pollak \cite{GP71} obtained the following result:
%%%%%%%%%%%%%%%%%%%%%%%%%%%%%%%%%%%%%%%%%%
%%%%%%%%%%%%% Theorem 1.1
%%%%%%%%%%%%%%%%%%%%%%%%%%%%%%%%%%%%%%%%%%
\begin{thm}[Graham and Pollak \cite{GP71}]
Let $T$ be a simple tree with $n$ vertices. Then
$$\det(D(T))=-(n-1)(-2)^{n-2},\eqno{(3)}$$
which is independent of the structure of $T$.
\end{thm}

Other proofs of Theorem 1.1 can be found in
\cite{Bapat051,Bapat052,BKN05,EGG76,GHH77,GL78,YY06}. In
particular, in \cite{YY06} we gave a simple method to prove $(3)$.
If $T$ is a weighted tree, Bapat, Kirkland, and Neumann
\cite{BKN05} generalized the result in Theorem 1.1 as follows.
%%%%%%%%%%%%%%%%%%%%%%%%%%%%%%%%%%%%%%%%%%
%%%%%%%%%%%%% Theorem 1.2
%%%%%%%%%%%%%%%%%%%%%%%%%%%%%%%%%%%%%%%%%%
\begin{thm}[Bapat, Kirkland, and Neumann \cite{BKN05}]
Let $T$ be a weighted tree with $n$ vertices and with edge weights
$\alpha_1,\alpha_2,\ldots,\alpha_{n-1}$. Then, for any real number
$x$,
$$\det(D(T)+xJ)=(-1)^{n-1}2^{n-2}\left(\prod_{i=1}^{n-1}\alpha_i\right)
\left(2x+\sum_{i=1}^{n-1}\alpha_i\right),\eqno{(4)}$$ where $J$ is
an $n\times n$ matrix with all entries equal to one.
\end{thm}

A direct consequence of Theorem 1.2 is the following:
%%%%%%%%%%%%%%%%%%%%%%%%%%%%%%%%%%%%%%%%%%
%%%%%%%%%%%%% Corollary 1.1
%%%%%%%%%%%%%%%%%%%%%%%%%%%%%%%%%%%%%%%%%%
\begin{cor}[Bapat, Kirkland, and Neumann \cite{BKN05}]
Let $D(T)$ be as in Theorem 1.2. Then
$$\det(D(T))=(-1)^{n-1}2^{n-2}\left(\prod_{i=1}^{n-1}\alpha_i\right)
\left(\sum_{i=1}^{n-1}\alpha_i\right).\eqno{(5)}$$
\end{cor}

Suppose $T$ is a weighted tree with the vertex set
$V(T)=\{v_1,v_2,\ldots,v_n\}$, and suppose the distance $d(u,v)$
between two vertices $u$ and $v$ is $\alpha$. Define two kinds of
$q$-{\sl distances} between $u$ and $v$, denoted by $d_q(u,v)$ and
$d^*_q(u,v)$, as $[\alpha]$ and $q^{\alpha}$ respectively, where
$$ [\alpha]=\left\{\begin{array}{ll}
\frac{1-q^{\alpha}}{1-q} &\ \ \mbox{if}\ \ q\neq 1,\\
\alpha &\ \ \mbox{otherwise}.
\end{array}
\right. $$ By definition, $[0]=0$ and
$[\alpha]=1+q+q^2+\ldots+q^{\alpha-1}$ if $\alpha$ is a positive
integer. We define two $q$-{\sl distance matrices} on the weighted
tree $T$, denoted by $D_q(T)$ and $D^*_q(T)$, as the $n\times n$
matrices with their $(i,j)$-entries equal to $d_q(v_i,v_j)$ and
$d^*_q(v_i,v_j)$, respectively. If $q=1$ then $D_q(T)$ is the
distance matrix $D(T)$ of $T$. Hence the distance matrix is a
special case of the $q$-distance matrix $D_q(T)$.

In quantum chemistry, if $T$ is a simple tree with vertex set
$V(T)=\{v_1,v_2,\ldots,v_n\}$,
$$W(T,q)=\sum\limits_{i<j}d_q^*(v_i,v_j)=\sum_{\{u,v\}\subseteq
V(T)}q^{d(u,v)}$$ is called the {\sl Wiener polynomial} of $T$
\cite{H71}, $D_1(T)$ is called the {\sl Wiener matrix}
\cite{GKYY06}, and the $q$-derivative $W'(T,1)$ is defined as the
{\sl Wiener index} of $T$ \cite{W471,W472}. The study of the
Wiener index, one of the molecular-graph-based structure
descriptors (so-called ``topological indices"), has been
undergoing rapid expansion in the last few years (see for example
\cite{KMPT92,RZ01,SYZ96,YY062,YY95}).

In the next section, we compute the determinants of $D^*_q(T)$ and
$D_q(T)$, and show that they are independent of  the structure of
$T$, and hence we generalize the results obtained by Graham and
Pollak \cite{GP71} and by Bapat, Kirkland, and Neumann
\cite{BKN05}. In Section 3, based on the results in Section 2, we
prove that the generating functions $F_n(q)=\sum\limits_{k\geq
0}N_{n,k}(T)q^k$ and $G_n(q)=\sum\limits_{k\geq 0}M_{n,k}(T)q^k$
of $\{N_{n,k}(T)\}_{k\geq 0}$ and $\{M_{n,k}(T)\}_{k\geq 0}$, as
defined in Problem 1.2, are exactly $\det(D_q^*(T))$ and
$\det(D_q(T))$, respectively. Hence, both $F_n(q)$ and $G_n(q)$
are independent of the structure of $T$, and this leads to a
resolution of Problem 1.2.
%%%%%%%%%%%%%%%%%%%%%%%%%%%%%%%%%%%%%%%%%%
%%%%%%%%%%%%% Section 2
%%%%%%%%%%%%%%%%%%%%%%%%%%%%%%%%%%%%%%%%%%
\section{Determinants of $D^*_q(T)$ and
$D_q(T)$} \hspace*{\parindent} First we compute the determinant of
$D_q^*(T)$.
%%%%%%%%%%%%%%%%%%%%%%%%%%%%%%%%%%%%%%%%%%
%%%%%%%%%%%%% Theorem 2.3
%%%%%%%%%%%%%%%%%%%%%%%%%%%%%%%%%%%%%%%%%%
\begin{thm}
Let $T$ be a weighted tree with $n$ vertices and with edge weights
$\alpha_1,\alpha_2,\ldots,\alpha_{n-1}$. Then, for any $n\geq 2$,
$$\det(D^*_q(T))=\prod_{i=1}^{n-1}(1-q^{2\alpha_i}), \eqno{(6)}
$$
which is independent of  the structure of $T$.
\end{thm}

{\bf Proof}\ \ We prove the theorem by induction on $n$. It is
trivial to show that the theorem holds for $n=2$ or $n=3$. Hence
we assume that $n\geq 4$. Without loss of generality, we suppose
that $v_1$ is a pendant vertex and $e=(v_1,v_s)$ is a pendant edge
with weight $\alpha_1$ in $T$. Let $d_i$ denote the $i$-th column
of $D_q^*(T)$ for $1\leq i\leq n$. Note that each entry along the
diagonal is one. Hence, by the definition of $D_q^*(T)$, we have
$$(d_1-q^{\alpha_1}d_s)^T=(1-q^{2\alpha_1},0,\ldots,0).$$
Thus
$$\det(D_q^*(T))=\det(d_1-q^{\alpha_1}d_s,d_2,d_3,\ldots,d_n)
=(1-q^{2\alpha_1})\det(D_q^*(T)^1_1),\eqno{(7)}$$ where
$D_q^*(T)^1_1$ equals $D_q^*(T-v_1)$. By induction, the theorem is
immediate from $(7)$. $\hfill\square$
%%%%%%%%%%%%%%%%%%%%%%%%%%%%%%%%%%%%%%%%%%
%%%%%%%%%%%%% Corollary 2.2
%%%%%%%%%%%%%%%%%%%%%%%%%%%%%%%%%%%%%%%%%%
\begin{cor}
Let $T$ be a simple tree with $n$ vertices. Then
$$\det(D^*_q(T))=(1-q^2)^{n-1},$$
which is independent of  the structure of $T$.
\end{cor}

To evaluate the determinant of $D_q(T)$ we must introduce some
terminology and notation. Let $A=(a_{ij})_{n\times n}$ be an
$n\times n$ matrix, and let $I=\{i_1,i_2,\ldots,i_l\}$ and
$J=\{j_1,j_2,\ldots,j_l\}$ be two subsets of $\{1,2,\ldots,n\}$.
We use $A^{i_1i_2\ldots i_l}_{j_1j_2\ldots j_l}$ to denote the
submatrix of $A$ by deleting rows in $I$ and columns in $J$.

Zeilberger \cite{Zeilb97} gave an elegant combinatorial proof of
Dodgson's determinant-evaluation rule \cite{Dodg1866} as follows:
$$\det(A)\det(A^{1n}_{1n})=\det(A_1^1)\det(A_n^n)-
\det(A_1^n)\det(A_n^1),\eqno{(8)}$$ where $A$ is a matrix of order
$n>2$. Let\\
$F(\alpha_1,\alpha_2,\ldots,\alpha_{n-1})$\\
$$=\frac{[\alpha_1][\alpha_{2}][\alpha_1+\alpha_{2}]}
{[2\alpha_1][2\alpha_{2}]}+\frac{[\alpha_{n-2}][\alpha_{n-1}][\alpha_{n-2}+\alpha_{n-1}]}
{[2\alpha_{n-2}][2\alpha_{n-1}]}+\sum\limits_{i=1}^{n-3}\frac{[\alpha_i][\alpha_{i+2}][\alpha_i+\alpha_{i+2}]}
{[2\alpha_i][2\alpha_{i+2}]}.$$

It is not difficult to prove the following lemma.
%%%%%%%%%%%%%%%%%%%%%%%%%%%%%%%%%%%%%%%%%%
%%%%%%%%%%%%% Lemma 2.1
%%%%%%%%%%%%%%%%%%%%%%%%%%%%%%%%%%%%%%%%%%
\begin{lem}
(a)\ If\ \ $n\geq 3$, $F(\alpha_1,\alpha_2,\ldots,\alpha_{n-1})$
is a symmetric function on
$\alpha_1,\alpha_2,\ldots,\alpha_{n-1}$.\\
(b)\ If\ \ $T$ is a weighted tree with two vertices and with edge
weight $\alpha_1$, $\det(D_q(T))=-[\alpha_1]^2.$\\
(c)\ If\ \ $T$ is a weighted tree with three vertices and with
edge weights $\alpha_1, \alpha_2$,
$\det(D_q(T))=2[\alpha_1][\alpha_2][\alpha_1+\alpha_2].$
\end{lem}
%%%%%%%%%%%%%%%%%%%%%%%%%%%%%%%%%%%%%%%%%%
%%%%%%%%%%%%% Theorem 2.4
%%%%%%%%%%%%%%%%%%%%%%%%%%%%%%%%%%%%%%%%%%
\begin{thm}
Let $T$ be a weighted tree with $n$ vertices and with edge weights
$\alpha_1,\alpha_2,\ldots,\alpha_{n-1}$. Then, for any $n\geq 4$,

$\det(D_q(T))=(-1)^{n-1}\left(\prod\limits_{i=1}^{n-1}[2\alpha_i]\right)\times$
$$
\left(\frac{[\alpha_1][\alpha_{2}][\alpha_1+\alpha_{2}]}
{[2\alpha_1][2\alpha_{2}]}+\frac{[\alpha_{n-2}][\alpha_{n-1}][\alpha_{n-2}+\alpha_{n-1}]}
{[2\alpha_{n-2}][2\alpha_{n-1}]}+\sum_{i=1}^{n-3}\frac{[\alpha_i][\alpha_{i+2}][\alpha_i+\alpha_{i+2}]}
{[2\alpha_i][2\alpha_{i+2}]}\right),\eqno{(9)}$$ which is
independent of the structure of $T$.
\end{thm}

{\bf Proof}\ \ We prove the theorem  by induction on $n$. Note
that there exist two trees with four vertices: the star $K_{1,3}$
and the path $P_4$. Let the edge weights of two weighted trees
$K_{1,3}$ and $P_4$ with four vertices be as shown in Figure 1 (a)
and (b), respectively. The $q$-distance matrices $D_q(T_{1,3})$
and $D_q(P_4)$ of $K_{1,3}$ and $P_4$ are as follows:
$$
D_q(T_{1,3})=\left(
\begin{array}{cccc}
\,0 & [\alpha_1+\alpha_2] & [\alpha_1+\alpha_3] & [\alpha_1] \\

[\alpha_1+\alpha_2] & 0 & [\alpha_2+\alpha_3] & [\alpha_2]\\

[\alpha_1+\alpha_3] & [\alpha_2+\alpha_3] & 0 & [\alpha_3]\\

[\alpha_1] & [\alpha_2] & [\alpha_3] & 0
\end{array}
\right)
$$
and
$$
D_q(P_4)=\left(
\begin{array}{cccc}
\,0 & \alpha_1 & [\alpha_1+\alpha_2] & [\alpha_1+\alpha_2+\alpha_3] \\

[\alpha_1] & 0 & [\alpha_2] & [\alpha_2+\alpha_3]\\

[\alpha_1+\alpha_2] & [\alpha_2] & 0 & [\alpha_3]\\

[\alpha_1+\alpha_2+\alpha_3] & [\alpha_2+\alpha_3] & [\alpha_3] &
0
\end{array}
\right).
$$
%%%%%%%%%%%%%%%%%%%%%%%%%%%%%%%%%%%%%%%%%%
%%%%%%%%%%%%% Figure 1
%%%%%%%%%%%%%%%%%%%%%%%%%%%%%%%%%%%%%%%%%%
\begin{figure}[htbp]
  \centering
 \scalebox{1}{\includegraphics{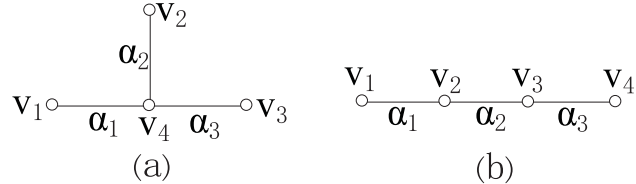}}
  \caption{\ (a)\ The weighted tree $T_{1,3}$.
  \ (b)\ The weighted tree $P_4$.}
\end{figure}

We calculate
$$\det(D_q(K_{1,3}))=\det(D_q(P_4))=$$$$-[2\alpha_1][2\alpha_2][2\alpha_3]
\left(\frac{[\alpha_1][\alpha_{2}][\alpha_1+\alpha_{2}]}
{[2\alpha_1][2\alpha_{2}]}+\frac{[\alpha_{2}][\alpha_{3}][\alpha_{2}+\alpha_{3}]}
{[2\alpha_{2}][2\alpha_{3}]}+\frac{[\alpha_1][\alpha_{3}][\alpha_1+\alpha_{3}]}
{[2\alpha_1][2\alpha_{3}]}\right).$$ Hence the theorem holds for
$n=4$.

Now we assume that $T$ is a weighted tree with $n$ vertices and
$n\geq 5$. We denote the $q$-distance matrix $D_q(T)$ of $T$ by
$D$. Note that $T$ has least two pendant vertices. Without loss of
generality, we assume both $v_1$ and $v_n$ are pendant vertices of
$T$. The unique neighbor of $v_1$ (resp. $v_n$) is denoted by
$v_s$ (resp. $v_{t}$). For convenience, we may suppose that the
weights of two edges $v_1v_s$ and $v_nv_t$ are $\beta_1$ and
$\beta_{n-1}$, and the weights of the edges in $T-v_1-v_n$ are
$\beta_2,\beta_3,\ldots,\beta_{n-2}$. Obviously,
$\{\beta_1,\beta_2,\ldots,\beta_{n-1}\}=
\{\alpha_1,\alpha_2,\ldots,\alpha_{n-1}\}$
($\{\beta_1,\beta_2,\ldots,\beta_{n-1}\}$ may be a multiset). Let
$d_i$ denote the $i$-th column of $D_q(T)$. By the definition of
$v_1$, $v_s$, $v_t$, and $v_n$, we have
$$(d_1-q^{\beta_1}d_s)^T=(-q^{\beta_1}[\beta_1],[\beta_1],
[\beta_1],\ldots,[\beta_1])$$ and
$$(d_n-q^{\beta_{n-1}}d_t)^T=([\beta_{n-1}],
[\beta_{n-1}],\ldots,[\beta_{n-1}],-q^{\beta_{n-1}}[\beta_{n-1}]),
$$ which imply the following:
$$\overline{d_1}^T=(d_1-q^{\beta_1}d_s)^T+\frac{-[\beta_1]}
{[\beta_{n-1}]}(d_n-q^{\beta_{n-1}}d_t)^T=
(-[2\beta_1],0,0,\ldots,0,(1+q^{\beta_{n-1}})[\beta_1]),$$ where
$d_1^T$ denotes the transpose of $d_1$. Hence
$$\det(D)=\det(d_1,d_2,\ldots,d_n)=
\det(\overline{d_1},d_2,d_3,\ldots,d_{n-1},d_n).$$ So we have
$$\det(D)=-[2\beta_1]\det(D_1^1)+(-1)^{n+1}
(1+q^{\beta_{n-1}})[\beta_1]\det(D_1^n).\eqno{(10)}$$ Similarly,
we have
$$\det(D)=-[2\beta_{n-1}]\det(D_n^n)+(-1)^{n+1}
(1+q^{\beta_{1}})[\beta_{n-1}]\det(D_n^1).\eqno{(10')}$$ On the
other hand, by Dodgson's determinant-evaluation rule $(6)$, we
have
$$\det(D)\det(D_{1n}^{1n})=\det(D_1^1)\det(D_n^n)-\det(D_1^n)\det(D_n^1). \eqno{(11)}$$
By the definition of the $q$-distance matrix $D$ (=$D_q(T)$) of
$T$, $\det(D_1^n)=\det(D_n^1)$. In particular, $D_1^1$, $D_n^n$,
and $D_{1n}^{1n}$ denote the $q$-distance matrices $D_q(T-v_1)$,
$D_q(T-v_n)$, and $D_q(T-v_1-v_n)$ of trees $T-v_1$, $T-v_n$, and
$T-v_1-v_n$, respectively. Note that $T-v_1$ (resp. $T-v_n$) is a
weighted tree with $n-1$ vertices and with edge weights
$\beta_2,\beta_3,\ldots,\beta_{n-1}$ (resp.
$\beta_1,\beta_2,\ldots,\beta_{n-2}$). Hence, by induction, we
have

$\det(D_1^1)=(-1)^{n-2}\left(\prod\limits_{i=2}^{n-1}[2\beta_i]\right)\times$
$$
\left(\frac{[\beta_2][\beta_{3}][\beta_2+\beta_{3}]}
{[2\beta_2][2\beta_{3}]}+\frac{[\beta_{n-2}][\beta_{n-1}][\beta_{n-2}+\beta_{n-1}]}
{[2\beta_{n-2}][2\beta_{n-1}]}+\sum_{i=2}^{n-3}\frac{[\beta_i][\beta_{i+2}][\beta_i+\beta_{i+2}]}
{[2\beta_i][2\beta_{i+2}]}\right) \eqno{(12)}$$ and

$\det(D_n^n)=(-1)^{n-2}\left(\prod\limits_{i=1}^{n-2}[2\beta_i]\right)\times$
$$
\left(\frac{[\beta_1][\beta_{2}][\beta_1+\beta_{2}]}
{[2\beta_1][2\beta_{2}]}+\frac{[\beta_{n-3}][\beta_{n-2}][\beta_{n-3}+\beta_{n-2}]}
{[2\beta_{n-3}][2\beta_{n-2}]}+\sum_{i=1}^{n-4}\frac{[\beta_i][\beta_{i+2}][\beta_i+\beta_{i+2}]}
{[2\beta_i][2\beta_{i+2}]}\right).\eqno{(13)}$$ Similarly,

$\det(D_{1n}^{1n})=(-1)^{n-3}\left(\prod\limits_{i=2}^{n-2}[2\beta_i]\right)\times$
$$
\left(\frac{[\beta_2][\beta_{3}][\beta_2+\beta_{3}]}
{[2\beta_2][2\beta_{3}]}+\frac{[\beta_{n-3}][\beta_{n-2}][\beta_{n-3}+\beta_{n-2}]}
{[2\beta_{n-3}][2\beta_{n-2}]}+\sum_{i=2}^{n-4}\frac{[\beta_i][\beta_{i+2}][\beta_i+\beta_{i+2}]}
{[2\beta_i][2\beta_{i+2}]}\right). \eqno{(14)}$$ From $(10)$ and
$(10')$,
$$[\det(D)]^2+[2\beta_1]\det(D)\det(D_1^1)+[2\beta_{n-1}]\det(D)\det(D^n_n)+$$$$
[2\beta_1][2\beta_{n-1}]\det(D^1_1)\det(D_n^n)=[2\beta_1][2\beta_{n-1}]\det(D^1_n)\det(D^n_1),$$
and hence by $(11)$ we have
$$[\det(D)]^2+[2\beta_1]\det(D)\det(D_1^1)+[2\beta_{n-1}]\det(D)\det(D^n_n)+
[2\beta_1][2\beta_{n-1}]\det(D)\det(D_{1n}^{1n})=0. \eqno{(12)}$$
Note that, by Theorem 1.1, if $q=1$ and $\beta_i=1$ for $1\leq
i\leq n-1$, then $\det(D)=-(n-1)(-2)^{n-1}$, which implies that
$\det(D)\neq 0$. Then by $(12)$ we have
$$\det(D)+[2\beta_1]\det(D_1^1)+[2\beta_{n-1}]\det(D^n_n)+
[2\beta_1][2\beta_{n-1}]\det(D_{1n}^{1n})=0. \eqno{(16)}$$

From $(12),(13),(14)$ and $(16)$, it is immediate that

$\det(D)=(-1)^{n-1}\left(\prod\limits_{i=1}^{n-1}[2\beta_i]\right)\times$
$$
\left(\frac{[\beta_1][\beta_{2}][\beta_1+\beta_{2}]}
{[2\beta_1][2\beta_{2}]}+\frac{[\beta_{n-2}][\beta_{n-1}][\beta_{n-2}+\beta_{n-1}]}
{[2\beta_{n-2}][2\beta_{n-1}]}+\sum_{i=1}^{n-3}\frac{[\beta_i][\beta_{i+2}][\beta_i+\beta_{i+2}]}
{[2\beta_i][2\beta_{i+2}]}\right).\eqno{(17)}$$

Note that
$\{\alpha_1,\alpha_2,\ldots,\alpha_{n-1}\}=\{\beta_1,\beta_2,\ldots,\beta_{n-1}\}$.
The theorem follows immediately from $(a)$ in Lemma 2.1 and
$(17)$. $\hfill\square$

Let $T$ be a weighted tree with the vertex set
$V(T)=\{v_1,v_2,\ldots,v_n\}$ and with the edge weights
$\alpha_1,\alpha_2,\ldots,\alpha_{n-1}$, and let $v_1$ and $v_n$
be two pendant vertices of $T$. The unique neighbor of $v_1$
(resp. $v_n$) is denoted by $v_s$ (resp. $v_{t}$). The proof above
also implies that
$$\det(D_q(T)_1^n)=[\alpha_1][\alpha_{n-1}]\prod_{i=2}^{n-2}[2\alpha_i],$$
where $\alpha_1$ and $\alpha_{n-1}$ are the weights of edges
$v_1v_s$ and $v_nv_t$, respectively.

If we set $q=1$ then the right hand side of $(9)$ in Theorem 2.4
equals
$$(-1)^{n-1}\prod_{i=1}^{n-1}(2\alpha_i)
\left(\frac{\alpha_1\alpha_2(\alpha_1+\alpha_2)}{(2\alpha_1)(2\alpha_2)}+
\frac{\alpha_{n-2}\alpha_{n-1}(\alpha_{n-2}+\alpha_{n-1})}{(2\alpha_{n-2})(2\alpha_{n-1})}
+\sum_{i=1}^{n-3}\frac{\alpha_i\alpha_{i+2}(\alpha_i+\alpha_{i+2})}{(2\alpha_i)(2\alpha_{i+2})}\right)$$
$$=(-1)^{n-1}2^{n-2}\left(\prod_{i=1}^{n-1}\alpha_i\right)
\left(\sum_{i=1}^{n-1}\alpha_i\right),$$ which implies Corollary
1.1 is a special case of Theorem 2.4.  Hence we generalize the
results obtained by Graham and Pollak \cite{GP71}, and by Bapat,
Kirkland, and Neumann \cite{BKN05}. In particular, the following
corollary is immediate from Theorem 2.4.

%%%%%%%%%%%%%%%%%%%%%%%%%%%%%%%%%%%%%%%%%%
%%%%%%%%%%%%% Corollary 2.3
%%%%%%%%%%%%%%%%%%%%%%%%%%%%%%%%%%%%%%%%%%
\begin{cor}
Let $T$ be a simple tree with $n$ vertices. Then
$$\det(D_q(T))=(-1)^{n-1}(n-1)(1+q)^{n-2},$$
which is independent of  the structure of $T$.
\end{cor}
%%%%%%%%%%%%%%%%%%%%%%%%%%%%%%%%%%%%%%%%%%
%%%%%%%%%%%%% Section 3
%%%%%%%%%%%%%%%%%%%%%%%%%%%%%%%%%%%%%%%%%%
\section{The quantities $M_{n,k}(T)$ and $N_{n,k}(T)$} \hspace*{\parindent}
Let $T$ be a simple tree and and $\mathcal A_{n,k}(T)=\{\sigma\in
\mathcal S_n |\ |\sigma_T|=k\}$. Partition $\mathcal S_n$ into
$\mathcal S_n=\mathcal A_{n,0}(T)\cup \mathcal A_{n,1}(T)\cup
\ldots \cup \mathcal A_{n,k}(T)\cup \ldots$.
%%%%%%%%%%%%%%%%%%%%%%%%%%%%%%%%%%%%%%%%%%
%%%%%%%%%%%%% Theorem 3.5
%%%%%%%%%%%%%%%%%%%%%%%%%%%%%%%%%%%%%%%%%%
\begin{thm}
Let $T$ be a simple tree with vertex set $\{v_1,v_2,\ldots,v_n\}$,
and let $N_{n,k}(T)$ be defined as in Problem 1.2. Then
$$
N_{n,k}(T)=\sum_{\sigma \in \mathcal
A_{n,k}(T)}\mathrm{sgn}(\sigma)=\left\{
\begin{array}{cc} 0 & \mbox{if}\ \ k \ \mbox{is \ odd,}\\
(-1)^{\frac{k}{2}}{n-1 \choose {\frac{k}{2}}} & \ \mbox{if}\ \ k\
\mbox{is\ even},
\end{array}
\right.$$ which is independent of  the structure of $T$.
\end{thm}
\par {\bf Proof}\ \ Let $F_n(q)=\sum\limits_{k\geq 0}N_{n,k}(T)q^k$ be
the generating function of $\{N_{n,k}(T)\}_{k\geq 0}$. Hence
$$
\begin{array}{lll}
F_n(q)&=&\sum\limits_{k\geq 0}\left (\sum\limits_{\sigma\in
\mathcal A_{n,k}(T)}\mathrm{sgn}(\sigma)\right)q^k
=\sum\limits_{k\geq 0}\left (\sum\limits_{\sigma\in \mathcal
A_{n,k}(T)}\mathrm{sgn}(\sigma)\right)q^{|\sigma_T|}\\
&=&\sum\limits_{k\geq 0}\left (\sum\limits_{\sigma\in \mathcal
A_{n,k}(T)}\mathrm{sgn}(\sigma)\right)q^{\sum\limits_{i=1}^nd(v_i,v_{\sigma(i)})}
=\sum\limits_{\sigma\in \mathcal S_n}\left (
\mathrm{sgn}(\sigma)q^{\sum\limits_{i=1}^nd(v_i,v_{\sigma(i)})}\right)\\
&=&\sum\limits_{\sigma\in \mathcal S_n}\left (
\mathrm{sgn}(\sigma)\prod\limits_{i=1}^nd_q^*(v_i,v_{\sigma(i)})\right).
\end{array}
$$
By the definition of $D_q^*(T)$, we have
$$\det(D_q^*(T))=\sum\limits_{\sigma\in \mathcal S_n}\left (
\mathrm{sgn}(\sigma)\prod\limits_{i=1}^nd_q^*(v_i,v_{\sigma(i)})\right).$$
The theorem is immediate from Corollary 2.2. $\hfill\square$

With notation as in the introduction, we state and prove our last
result.
%%%%%%%%%%%%%%%%%%%%%%%%%%%%%%%%%%%%%%%%%%
%%%%%%%%%%%%% Theorem 3.6
%%%%%%%%%%%%%%%%%%%%%%%%%%%%%%%%%%%%%%%%%%
\begin{thm}
Let $T$ be a simple tree with vertex set $\{v_1,v_2,\ldots,v_n\}$,
and let $M_{n,k}(T)$ and $\phi_{\sigma,k}(T)$ be as in (2). Then
$$
M_{n,k}(T)=\sum_{\sigma \in \mathcal
S_{n}}\mathrm{sgn}(\sigma)\phi_{\sigma,k}(T)=(-1)^{n-1}(n-1){n-2\choose
k},$$ which is independent of the structure of $T$.
\end{thm}

{\bf Proof}\ \ Let $G_n(q)=\sum\limits_{k\geq 0}M_{n,k}(T)q^k$ be
the generating function of $\{M_{n,k}(T)\}_{k\geq 0}$. Hence
$$
\begin{array}{lll}
G_n(q)&=&\sum\limits_{k\geq 0}\left (\sum\limits_{\sigma\in
\mathcal S_n}\mathrm{sgn}(\sigma)\phi_{\sigma,k}(T)\right)q^k
=\sum\limits_{\sigma\in \mathcal
S_n}\left(\mathrm{sgn}(\sigma)\sum\limits_{k\geq
0}\phi_{\sigma,k}(T)q^k\right)\\
&=&\sum\limits_{\sigma\in \mathcal
S_n}\mathrm{sgn}(\sigma)(1+q+\ldots+q^{d(v_1,v_{\sigma(1)})-1})
\ldots(1+q+\ldots+q^{d(v_n,v_{\sigma(n)})-1})\\
&=&\sum\limits_{\sigma\in \mathcal
S_n}\mathrm{sgn}(\sigma)d_q(v_1,v_{\sigma(1)})
d_q(v_2,v_{\sigma(2)})\ldots d_q(v_n,v_{\sigma(n)}).
\end{array}
$$
By the definition of $D_q(T)$, we have
$$\det(D_q(T))=\sum\limits_{\sigma\in \mathcal
S_n}\mathrm{sgn}(\sigma)d_q(v_1,v_{\sigma(1)})
d_q(v_2,v_{\sigma(2)})\ldots d_q(v_n,v_{\sigma(n)}).$$ The theorem
follows immediately from Corollary 2.3. $\hfill\square$

By Remark 1.1, Theorems 3.5, and 3.6,
$M_{n,k}=(-1)^{n-1}(n-1){n-2 \choose k}$, while $N_{n,k}=0$ if $k$
is odd and $N_{n,k}=(-1)^{\frac{k}{2}}{n-1\choose {\frac{k}{2}}}$
otherwise.

Our method to prove Theorems 3.5 and 3.6 is completely algebraic.
Therefore it would be interesting to consider the following
problem.
%%%%%%%%%%%%%%%%%%%%%%%%%%%%%%%%%%%%%%%%%%
%%%%%%%%%%%%% Problem
%%%%%%%%%%%%%%%%%%%%%%%%%%%%%%%%%%%%%%%%%%
\begin{prob}
Give combinatorial proofs of Theorems 3.5 and 3.6.
\end{prob}
%%%%%%%%%%%%%%%%%%%%%%%%%%%%%%%%%%%%%%%%%%%
%%%%%%%%%%%%%%% Acknowledgements
%%%%%%%%%%%%%%%%%%%%%%%%%%%%%%%%%%%%%%%%%%%
\vskip0.5cm \noindent {\bf Acknowledgements}
\par Thanks to Professor David B. Chandler and the referee for providing
many very helpful suggestions for revising this paper.
%\newpage
%%%%%%%%%%%%%%%%%%%%%%%%%%%%%%%%%%%%%%%%%%%
%%%%%%%%%%%%%%% References
%%%%%%%%%%%%%%%%%%%%%%%%%%%%%%%%%%%%%%%%%%%

\end{document}